\documentclass{birkjour}%

\usepackage{graphicx}
\usepackage{latexsym}
\usepackage{amssymb,amsmath}
\usepackage{url}
\usepackage{multicol}
\usepackage{color} % Temporary, for making annotations and comments.

 \theoremstyle{definition}
 
 \theoremstyle{remark}

 \numberwithin{equation}{section}

\usepackage{txfonts}%%!!
\makeatletter%
\input{ot1txtt.fd}
\makeatother%

\newcommand{\R}{\mathbb{R}}	
\newcommand{\C}{\mathbb{C}}
\newcommand{\HQ}{\mathbb{H}}

\newcommand{\be}{\begin{equation}}
\newcommand{\ee}{\end{equation}}

\renewcommand{\O}{\mathrm{O}}

\renewcommand{\eqref}[1]{(\ref{#1})}
\newcommand{\vect}[1]{\ensuremath{\mbox{\textbf{\textit{#1}}}}}
  \newcommand{\qi}{\ensuremath{\mbox{\boldmath $i$}}}
  
  \newcommand{\qj}{\ensuremath{\mbox{\boldmath $j$}}}
  
  \newcommand{\qk}{\ensuremath{\mbox{\boldmath $k$}}}

\newcommand{\change}[1]{{#1}}
\newcommand{\boldx}{\mbox{\boldmath $x$}}
\newcommand{\changeNitta}[1]{{{#1}}}
\newcommand{\changeNittaRevise}[1]{{#1}}	% for revision

\begin{document}

\title{Applications of Clifford's Geometric Algebra}

%----------Author 1
\author{Eckhard Hitzer}
\address{%
College of Liberal Arts, Department of Material Science,\\ 
International Christian University,\\
Tokyo 181-8585,\\ 
Japan}
\email{hitzer@icu.ac.jp}
%\thanks{This work was completed with the support of our \TeX-pert.}
%----------Author 2
\author{Tohru Nitta}
\address{Mathematical Neuroinformatics Group,\\
    Human Technology Research Institute,\\
    National Institute of Advanced Industrial Science and Technology,\\
    AIST Tsukuba Central 2,\\
    Ibaraki 305-8568,\\
    Japan}
\email{tohru-nitta@aist.go.jp}
%----------Author 3
\author{Yasuaki Kuroe}

\address{%
Department of Information Science,\\
Kyoto Institute of Technology,\\ 
Matsugasaki, Sakyo-ku,\\ 
Kyoto 606-8585,\\
Japan}
\email{kuroe@kit.ac.jp}
\subjclass{Primary 15A66; \\Secondary 68T40,62M45,68U10,60G35}
% 15A66   	Clifford algebras, spinors
% 68T40   	Robotics
% 62M45   	Neural nets and related approaches
% 68U10   	Image processing
% 60G35   	Signal detection and filtering
\keywords{Hypercomplex algebra, hypercomplex analysis, geometry, science, engineering, {applications}}

\date{March 23, 2012}
\dedicatory{Soli Deo Gloria}
%%% ----------------------------------------------------------------------

\begin{abstract}
We survey the development of Clifford's geometric algebra and some of its engineering applications during the last 15 years. Several recently developed applications and their merits are discussed in some detail. We thus hope to clearly demonstrate the benefit of developing problem solutions in a unified framework for algebra and geometry with the widest possible scope: from quantum computing and electromagnetism to satellite navigation, from neural computing to camera geometry, image processing, robotics and beyond.
\end{abstract}

%%% ----------------------------------------------------------------------
\maketitle
%%% ----------------------------------------------------------------------

\section{{Introduction}}

{
In this review we try to give an overview on how applications of Clifford's geometric algebra mainly in the areas of neural computing, image and signal processing, computer and robot vision, control problems and other areas have developed over the past 15 years. The enormous range of applications developed in the past decades makes a complete overview next to impossible. We therefore restricted the review to the above mentioned application fields. We further particularly focus on the representative proceedings of the following conferences: 
Applied Clifford Algebra in Cybernetics, Robotics, Image Processing and Engineering (ACACSE 1999) \cite{BCS:GAappScEng2001}, 
Applied Geometric Algebras in Computer Science and Engineering (AGACSE 2001) \cite{DDL:AGACSE2002}, AGACSE 2008 \cite{BCS:GACEngCS2010}, and AGACSE 2010 \cite{DL:GGAP2011}. 

But for giving a more complete picture, we added a number of further related references. This review is mostly descriptive. In order to provide deeper insights, we additionally chose to review several selected publications in more detail. Here our choice is necessarily somewhat subjective, but we hope that the readers can thus get a better idea of the general characteristics of solutions for applied science and technical problems with Clifford's geometric algebra. One such characteristic is, e.g., that due to the general coordinate invariance of the methods applied, the choice of coordinate system can often beneficially be postponed until concrete application data are to be processed. 

The review is organized as follows. Without any claim for completeness,
Section \ref{sc:MAlgCGA} briefly recounts how motor algebra and conformal geometric algebra arose.
Section \ref{sc:CliffNN} focuses on Clifford algebra neural computing, both with neural networks and with support vector machines. 
Section \ref{sc:SIProc} reviews the wide area of applications to signal and image processing, including electromagnetic signals, image geometry and structure, vector fields, and color image processing. 
This is continued with an overview of the applications of Clifford integral transformations of Fourier type, and other types like wavelets. Section \ref{sc:CompRobV} treats applications to computer and robot vision, i.e. for orientation, pose, motion and tracking, as well as for camera geometries and scene analysis. 
Section \ref{sc:KinDynRob} deals with applications to (mainly) kinematics and dynamics of robots. 
Section \ref{sc:CtrlProbl} shows some applications to control problems, in particular optical positioning and geometric constraint solving. 
Finally, Section \ref{sc:othApp} on other applications, contains an application to automatic theorem proving, which also leads to interesting practical results, in passing mentions aspects of the Clifford algebra treatment of quantum computing and light, and remarks (not exhaustively) on some of the software available for Clifford algebra. 
}

\section{Emergence of motor algebra (= dual quaternions) and conformal geometric algebra \label{sc:MAlgCGA}}

Clifford's geometric algebras $Cl(p,q,r)$ constitute an extension of real, complex and quaternion algebras to complete associative algebras of subspaces of vector spaces $\R^{p,q,r}$ \cite{EH:SICEintro}. Following established notation, we denote $\R^{p,q} := \R^{p,q,0}$, and $Cl(p,q):= Cl(p,q,0)$. Regarding examples, the
complex numbers are isomorphic to $Cl(0,1)$, $Cl^+(2,0)$, and $Cl^+(0,2)$, etc., and quaternions are isomorphic to $Cl(0,2)$, and $Cl^+(3,0)$, and $Cl^+(0,3)$, etc.

We will see that in general, apart from special applications (like visibility, where Grassmann algebra is sufficient) nowadays \textit{conformal geometric algebra} 
\linebreak
$Cl(4,1)$, abbreviated as CGA, or more general $Cl(p+1,q+1)$, has become the most widely applied system. Complex numbers, quaternions, (complex and dual) biquaternions and motor algebra are all included as subsystems (subalgebras) in conformal geometric algebra.

Note that the \textit{motor algebra} $Cl^+(3,0,1)$ is isomorphic to $Cl^+(0,3,1)$, and to Clifford's original 19th century biquaternions (dual quaternions). 

In \cite{BCS:GAappScEng2001} (the AGACSE 1999 proceedings) has four papers (chp. 1 by D. Hestenes, chp. 2 by G. Sobczyk, chp. 3 by J. M. Pozo and G. Sobczyk, chp. 4 by H. Li) introducing the conformal model in 
$Cl(p+1,q+1)$. Specialized to $Cl(n+1,1)$ it provides a universal model for conformal geometries of Euclidean, spherical and double-hyperbolic spaces. 
It establishes an algebraic formulation of conformal geometry according to the $(n+2)$D $(p,q)$-Horosphere model, $n=p+q$. It is also called double projective. Moreover, it constitutes an embedding of the previously known motor algebra $Cl^+(3,0,1)$. The most comprehensive exposition of conformal geometric algebra is perhaps given in the first four chapters of \cite{GS:GCompwGA2001} by D. Hestenes, H. Li, and A. Rockwood.

\change{%
As shown, e.g., by C. Gunn in chapter 15 of \cite{DL:GGAP2011}, also the homogeneous (or projective) algebra $Cl(p,0,1)$ can be favorably applied to problems in Euclidean geometry. This has the computational advantage, that only one extra dimension is needed instead of the two for conformal geometric algebra. A slight drawback is that in homogeneous algebra $Cl(p,0,1)$ transformations like rotations and translations do not apply to all geometric objects in the same single way that their versor representations do in conformal geometric algebra \cite{DFM:GAfCS}. We moreover note, that the algebra $Cl(p,0,1)$ also appears as a subalgebra of $Cl(p+1,1)$, by selecting only one null-vector from the Minkowski plane added to $Cl(p,0)$, and dropping the other. 
}

\change{%
For an in depth mathematical treatment of all hypercomplex algebras concerned in this review we kindly refer the interested reader to the comprehensive modern textbook by H. Li in \cite{HL:IAandGR}. We also note that for a deeper understanding, the reader should always examine the detailed algebraic definitions for the particular use of each hypercomplex algebra in every application reference, which we discuss here, because due to the overview character of this review, we cannot discuss the mathematical foundations in full detail. 
}

%%-----------------------------------------------------
%% 3. Clifford neural computing (rewritten by Tohru NittaNitta)
%%-----------------------------------------------------
\section{Clifford neural computing\label{sc:CliffNN}}

%--------------
% Introduction
%--------------
\changeNittaRevise{
\change{An a}rtificial neural network is a mathematical model inspired by biological neural networks, 
and has been recognized as a useful tool in machine learning. 
For example, neural networks can learn nonlinear relationships just using given input 
and output data, and can also be used as an associative memory. 
Neural networks are composed of a number of nonlinear computational elements (called neuron\change{s}) 
which operate in parallel and are arranged in a manner reminiscent of biological neural interconnections. 
Neural networks can be grouped into two categories: 
one is feedforward neural networks, in which each neuron has no loops, 
and the other is recurrent neural networks, where loops occur with feedback connections.}

This {section} describes several topics on neural networks extended to Clifford number called 
{\it Clifford neural networks} (Clifford NN). 
The study of extending neural networks to higher-dimensions such as complex numbers 
has remained active since the 1990s \cite{hirose2003,hirose2006,nit08a,nit09}. 

There appear to be two approaches for extending the usual neural networks to higher dimensions. 
One approach is to extend the number field, i.e., from real numbers $x$ (1 dimension), 
to \changeNitta{$2^nD$ multivectors\footnote{
Octonions by themselves are no Clifford algebra. But it is still possible to elegantly work with octonions in Clifford algebra, by embedding octonions in $Cl(0,7)$ \cite[Chapter 7.4]{PL:CAandS}.
} 
such as complex numbers $z=x+iy$ ($2^1$ dimensions) and quaternions 
{$q=a+\qi b+\qj c+\qk d$} ($2^2$ dimensions) (called {\it Clifford approach} here). }
% to {octonions} (8 dimensions), to sedenions (16 dimensions), $\cdots$. 
%
Another approach is just to extend the dimensionality of the weights and threshold values 
from 1 dimension to $n$ dimensions using $n$-dimensional real valued vectors. 
Moreover, the latter approach has two varieties: 
(a) weights are $n$-dimensional matrices \cite{nit92,nit94,nit06,nit07}, 
(b) weights are $n$-dimensional vectors \cite{koba04,nit93}. 
In this \changeNitta{section} we deal with the \changeNitta{Clifford approach}. 

\subsection{Clifford neural networks}

Studies of the Clifford NNs has been focused on multi-layered type ones. 
Pearson first proposes a Clifford neural network model \cite{pearson92,pearson94}, 
in which he formulates a multi-layered Clifford NN model, derived a Clifford back-propagation 
learning algorithm, and clarifies the ability on approximate functions. 

Buchholz formulates in his PhD thesis \cite{SB:PhD} another multi-layered Clifford NN 
in $Cl(p,q,r)$ with a split-type activation function which is different from that Pearson adopted. 
He shows that single neurons already describe geometric transformations, 
and studies the spinor Clifford neurons (SCN), whose weights act like rotors from two sides. 
Next Buchholz discusses Clifford multilayer perceptrons, and back-propagation learning algorithm, 
including approximation and prediction applications, and comparisons with conventional neural networks. 
\changeNittaRevise{
As a result, it is learned that the performance of the Clifford multilayer neural networks is superior to 
that of the usual real-valued neural networks.}
Important results of Buchholz' PhD thesis are summarized in \cite{BS:CNMLP}.

In part II of \cite{LOS:CAGAapp2005} G. Sommer gives an overview of applications 
of geometric algebra in robot vision (pp. 258--277) including the so called spinor Clifford 
neuron with rotor weights in $Cl^+(p,q)$, learning of Moebius transformations in $Cl(p+1,q+1)$, 
and the hypersphere neuron in $Cl(n+1,1)$. 
\changeNittaRevise{
\change{Note that} Moebius transformations cannot be learned with the usual real-valued neural networks.} 

J. Rivera-Rovelo, E. Bayro-Corrochano and R. Dillmann develop with conformal geometric 
algebra $Cl(4,1)$ in \cite{BCS:GACEngCS2010}, pp. 191--210, geometric neural computing 
for 2D contour and 3D surface reconstruction. B. Cruz, R. Barron and H. Sossa study 
in \cite{BCS:GACEngCS2010}, pp. 211--230, geometric associative memories with 
applications to pattern classification in $Cl(4,1)$. M. T. Pham et al discuss 
in \cite{BCS:GACEngCS2010}, pp. 231--248, the classification and clustering of spatial 
patterns with geometric algebra $Cl(n,0)$.

%------------------------
% Recurrent Clifford NN
%------------------------
Y. Kuroe first proposes models of recurrent Clifford NNs and discusses their
dynamics from the point of view of the existence of energy functions \cite{kuroe_ijcnn2011,kuroe_iconip2011}. 

%%-------------------------------------
%% 3.2 Clifford Support Vector Machine
%%-------------------------------------
\subsection{Clifford support vector machine}
\changeNittaRevise{
\change{A} Support Vector Machine (SVM) is a supervised learning system that can be applied to 
classification and regression problems, which has been proposed in 1995 \cite{Vapnik}. 
One of \change{the superior} points of the SVM is its generalization ability. 
\change{It} has been found \change{in} computer simulations \change{that SVMs} extended \change{to Clifford} algebra 
\change{reveal excellent} performance compared \change{to real}-valued SVMs, 
as described in the following literatures. 
}

In \cite{BCS:GAappScEng2001} (chp. 13 by E. Bayro-Corrochano and R. Vallejo) 
geometric feedforward neural networks and support vector machines are considered 
in $Cl(p,q,r)$ with examples in $Cl(2,0,0)$, $Cl(0,2,0)$, $Cl(3,0,0)$ and $Cl(3,0,1)$, 
see also \cite{EBC:GCPAS2001} by E. Bayro-Corrochano.

In \cite{BC:GNC} Support Multivector 
Machines (SMVM or Clifford SVM = CSVM) are introduced with a variety of kernels, 
able to construct optimally separating hyperplanes in multivector feature spaces. 
Then multivector regression is briefly considered (see also \cite{BC:CSVM}). 
The experimental part presents learning of nonlinear mappings, encoder-decoder tasks, 
3D pose estimation with 2D and 3D data, curvature encoding, and 3D rigid motion estimation. 

Additionally, in \cite{BC:CSVM} recurrent CSVM are studied, and CSVMs are evaluated 
for classification, regression and recurrence. Evaluation problems are 3D nonlinear 
geometric classification, 3D object class recognition including real robot data, 
multicase interpolation, time series data, and robotic navigation in discrete labyrinths.

\section{Applications to signal and image processing\label{sc:SIProc}}

\subsection{Electromagnetic signals}

\change{%
Regarding electromagnetism, spacetime algebra (STA) \cite{DH:STA} in $Cl(1,3)$ provides a unifying framework for all four Maxwell equations in the form of
\be
  \label{eq:Maxwell}
  \nabla F = J,
\ee
where $\nabla$ is the spacetime vector derivative, $F=E + i_4 B$ is the electromagnetic field bivector, $i_4=e_0e_1e_2e_3$ is the spacetime hypervolume pseudoscalar, and $J$ is the spacetime vector current. Decomposition into the various grade parts of \eqref{eq:Maxwell} gives all four Maxwell equations. The vector derivative operator $\nabla$ is invertible via integration and Greens functions, thus allowing to compute the fields in terms of the currents. 
}

\change{%
This new STA approach is used} in \cite{CgeoAPME1996} (chps. 8--11 by S. Gull, C. Doran, A. Lasenby) to treat electromagnetic waves, including boundaries, propagation in layered media, and tunneling. In \change{the alternative Algebra of Physical Space (APS) approach in} $Cl(3,0)$ parabivectors (scalars + bivectors) are used for polarized electromagnetic waves \cite{CgeoAPME1996} (chps. 17+18 by W. Baylis). 

\change{Complex vector function analysis in higher dimensions is a well established field. Historically a former student of D. Hilbert \cite{GHS:HolFctnD}, the Swiss mathematician R. Fueter \cite{RF:AnThQuatVar,RF:FunktTh} began to generalize the Cauchy-Riemann equations with the help of quaternions to higher dimensions, which has been subsequently extended to Clifford algebras of $n$-dimensional vector spaces \cite{GHS:HolFctnD} and leads to a mathematical refinement \cite{BKS:GenMDHilbCA} compared to the theory of differential forms. This new form of hypercomplex Clifford analysis provides new tools for applied problems in higher dimensions. Thus in}
\cite{GSp:QCCPE1997}, the authors K. G\"{u}rlebeck and W. Spr\"{o}ssig, develop based on quaternions, complex quaternions, and Clifford algebras $Cl(p,q)$, Clifford valued functions and forms, as well as a Clifford operator calculus including integral transformations like the Clifford Fourier transform. Next they consider boundary value problems for Clifford valued functions relevant for electromagnetism, thermo-elasticity, fluid dynamics and transmission problems. Finally they introduce discretization for numerical applications.

\subsection{Image geometry and structure}

\change{%
In image processing many algorithms apply differential invariants from Euclidean space to the "image surface" (in "image space", that is the picture plane times intensity space). In chp. 28 of \cite{DDL:AGACSE2002} J. Koenderink notes, that this makes little sense since these invariants are with respect to Euclidean isometries, e.g. rotations. But obviously it is impossible to rotate the image surface to see its other side, or the intensity to a spatial direction. Thus the angle measure cannot be periodic in planes other than the picture plane. There is a clear need to set up the proper transformation group to arrive at a set of meaningful invariants. Koenderink therefore proceeds to develop in $Cl(2,0,1)$ a generic framework for image geometry,} compare also chp. 36 of \cite{RA:CAappsMPE2004} by the same author. Koenderink represents similarities (=similarity transformations) of images by rotor transformations in $Cl(2,0,1)$ and researches differential image geometry of 1D curves and 2D surfaces in the 3D image space of the picture plane and an additional perpendicular log intensity axis, where this axis is represented by a null vector, squaring to 0.

% commented out by T.Nitta 2012.4.10
%A general computer generated 2D vector field model in $Cl(2,0)$ is presented by A. Rockwood and S. Binderwala chp. 16 of \cite{DDL:AGACSE2002}. 

\change{
In the case of multidimensional signals, quadrature filters can only be applied to a preferential direction with a need for prior orientation sampling and the application of steerable or orientable filters. 
In \cite{SZ:PerAct2000} (pp. 175--185) M. Felsberg and G. Sommer for the first time introduce a linear approach, based on geometric algebra, to obtain an isotropic analytic signal, with an amplitude independent of the local orientation. That means they} consider multidimensional isotropic generalizations of quadratic filters in geometric algebra $Cl(0,n)$, and test them on 2D and 3D signals. 
Moreover in chp. 38 of \cite{RA:CAappsMPE2004} M. Felsberg and G. Sommer derive a general operator for analyzing the local structure of 2D images, the structure multivector, by embedding the images in $Cl(3,0)$.

\change{Traditional ray tracing algorithms employ linear algebra matrix and vector algorithms in three dimensions, or four-dimensional projective geometry.}
In \cite{DF:M3DGeo} the performance of geometric algebra ray tracing algorithms in Euclidean $Cl(3,0)$, projective $Cl(4,0)$ and conformal $Cl(4,1)$ is compared with conventional linear algebra algorithms in 3D and in projective 4D space.

%%-------------------------------------------------
%% 4.3 Vector field ... created by T.Nitta 2012.4.12
%%-------------------------------------------------
\subsection{Vector field}

Vector {fields have} been often used {in computer} graphics and {for simulations}. 
{A} vector {field} in {a} 2-dimensional Euclid space $\R^2$ is defined as a smooth map: 
\begin{equation}
  \label{VF}
  E: \R^2 \rightarrow \R^2, \,\, \boldx = x e_1 + y e_2 \mapsto E_1 e_1 + E_2 e_2\, ,
\end{equation}
where $x, y, E_1, E_2\in \R$, and $e_1, e_2$ are basis vectors. 
A {critical} point of Eq. (\ref{VF}) occurs where $E(\boldx)=0$, 
which is used to characterize the vector field. 
{The} Poincar\'{e}-index or winding number is an invariant for critical points, 
and counts the number of turns of the field around the critical point. 
%%--------------------------
%%　Polya's method
%%--------------------------
Polya proposed an approach in which critical points are defined by roots and poles 
of a complex function \cite{braden1985,braden1991}. 
\begin{equation}
  E(z) 
  = \alpha \frac{\prod_i(z - p_i)^{r(i)}}{\prod_i(z - q_i)^{s(i)}}\, ,
\end{equation}
where $\alpha$ is a complex constant. 
{The} Poincar\'{e}-index of a critical point $p_i$ is defined as $r(i) \geq 0$, 
and that of a critical point $q_i$ as $-s(i) \leq 0$. 
In this connection, $p_i$ {is zero} of the order $r(i)$ of the complex function $E$, and $q_i$ is a pole 
of the order $s(i)$. 
%-------------------------------
%　Problems in the Polya's method
%-------------------------------
{However, Polya's} method has two problems. One is that numerical problems arise 
when calculating integral curves (a curve between critical points) near the pole. 
{Another} one is that a pole cannot {correctly} model a critical point with a negative index. 
%%------------------------
%% Conjugate field method
%%------------------------
Scheuermann et al. propose the conjugate field method that solves these problems by using 
Clifford algebras \cite{sch}. 
{The}
Clifford algebra that they use is $Cl(2,0)$ (four dimensions) whose basis are $1, e_1, e_2, i$, 
and satisfy the following rules: $1^2 = e_1^2 = e_2^2 = 1, i^2 = -1, e_1e_2 = - e_2e_1 = i$. 
The basis {element} $i$ has similar nature to the imaginary number of complex numbers, 
{therefore} ${z=} a+bi \in {Cl(2,0)}$ can be interpreted as {a} `complex number'. 
Let $z = x + 0 e_1 + 0 e_2 + yi \in Cl(2,0)$ for any $\boldx = x e_1 + y e_2 \in \R^2$. 
Then, we obtain $E(\boldx) = E(z, \bar{z}) e_1$ using the relation $e_2 = - i e_1$. 
In general, $E(z, \bar{z}) = \alpha F_1 \cdots F_n$ where $\alpha$ is a complex constant, and 
%
% \newline
% \comment{In the above paragraph without explanation two notations are used
% for basis vectors: $e_1,e_2$ 
% and $\sigma_1, \sigma_2$. If possible, this should be reduced to only one
% pair of basis vectors, e.g. only $e_1,e_2$. What do you think?}
%
% \comment{\bf I agree. I unified with $e_1$ and $e_2$. }
%
\begin{equation}
  F_i =
  \left\{
    \begin{array}{ll}
     (\bar{z} - (x_i - iy_i))^{r(i)} \hspace{0.8cm} \left(r(i) > 0 
    \right)\\
     (z - (x_i + iy_i))^{-r(i)} \hspace{0.6cm} \left(r(i) < 0 \right).
    \end{array}
  \right.
\end{equation}
$F_i$ expresses the critical point $(x_i, y_i)$ {with}  Poincar\'{e}-index $r(i)$. 
%%------------------------
%% Comparison
%%------------------------
{A} critical point with a negative Poincar\'{e}-index is expressed using poles {in} Polya's method 
while it is expressed as the root of the variable $z$ in the conjugate field method. 
This difference causes the stability and high speed of the calculations of the conjugate field method 
compared {to} Polya's method. 
%%------------------------------------
%% Rockwood's work (software)
%%------------------------------------
Rockwood et al. provide {a} computer program that shows integral curves and vector fields on the display 
using the conjugate field method, and can be modified interactively \cite{rock}.

\subsection{Color image processing}

\change{%
Conventional color image processing relies on marginal channel wise red, green and blue color processing. Thus color space is not approached in a holistic three-dimensional approach. New holistic multidimensional non-marginal treatments in quaternions and geometric algebra overcome this limitation and lead to a genuine color space image processing, which allows, e.g., to compare pictures rotated relative to each other in color space, to treat color illumination changes, etc.}

In \cite{SZ:PerAct2000} (pp. 78--103) V. Labunets, E. Labunets-Rundblad, and J. Astola develop a method to calculate invariants of $k$-multicolor $n$D images using multiplet Fourier-Clifford-Galois number theoretical transforms, thereby reducing computational complexity. 
In chp. 7 of \cite{GS:GCompwGA2001} E. Rundblad-Labunets and V. Labunets use Clifford algebras $Cl(p,q,r)$ as spatial-color Clifford algebras for invariant pattern image recognition. 
Color edge detection using rotors in $Cl(3,0)$ is discussed in chp. 29 of \cite{DDL:AGACSE2002} by E. Bayro-Corrochano and S. Flores. 
T. Batard, M. Berthier and C. Saint-Jean develop in \cite{BCS:GACEngCS2010}, pp. 135--162, a special color version of the Clifford Fourier transform in $Cl(4,0)$ for 2D color image processing. 

In \cite{MSJM:CObRCFT} techniques for 
color object recognition based on a special type of Clifford Fourier transform \cite{BBS:colCFT} are developed. The color CFT embeds a 3D color vector signal in a 4D space $\R^4$ and uses a bivector $B \in Cl(4,0)$ and its dual $I_4B$. It needs to be distinguished from the standard CFT introduced by \cite{ES:CFT} and \cite{HM:CFT}, which applies to general multivector-valued signals and has a well-defined general inverse transform. The color CFT of \cite{BBS:colCFT} lacks a general definition for multivector signals and an inverse. 

The first technique of \cite{MSJM:CObRCFT} adapts the generalized Fourier descriptors (GFD) of \cite{FSetal:GFD} for use with the color CFT (GCFD). The work claims, that the GCFDs are more compact, less complex and have better recognition rates (which strongly depends on the color image database used in experiments), than the GFDs of \cite{FSetal:GFD} with separate processing of the red, green and blue channels. Only two fast Fourier transforms are necessary compared to three in \cite{FSetal:GFD}. But it must be noted, that the authors use twice the number of descriptors for the GCFDs, than for the GFDs in their comparison.  

The second technique of \cite{MSJM:CObRCFT} introduces a color phase correlation with the same color CFT, based on a color CFT shift theorem. Both techniques are applied to several color image databases. The influence of the choice of the analysis bivector plane $B$ in color space is discussed. The work itself admits, that a mathematically sound phase correlation for color images needs an inverse CFT, which is still not available for the color CFT of \cite{BBS:colCFT}. 

In the context of phase correlation it might therefore be worthwhile to treat a color vector signal with the standard CFT of \cite{ES:CFT} and \cite{HM:CFT}, which is well-defined for general multivector-valued signals and has a general inverse transform. It should also be very promising to use of the generalized quaternion Fourier transform (QFT), called OPS-QFT (OPS = orthogonal 2D planes split of quaternions) \cite{HS:ICCA9, EH:ICNAAM2011}, which is well-defined for all quaternion signals, has a general inverse, is dimensionally minimal, embeds colors, and also has the feature of steerable 2D analysis planes.

\subsection{Clifford integral transformations}

\subsubsection{Clifford Fourier type transformations}

{\change{%
Clifford Fourier type integral transformations are an active field of research. Broadly speaking there currently are mainly three different, but related approaches. The first eigenfunction approach in complex Clifford algebras leads to new hypercomplex transformations by assuming certain eigenvalues for a function basis of a space of square integrable Clifford algebra valued functions \cite{BX:OnCFT,BF:VbVFTinCA}. The second approach replaces the imaginary unit $i \in \C$ by a Clifford algebra square root of $-1$ \cite{SJS:BiqR-1,HA:GR-1,HHA:SqR-1inrCA} taken from a real Clifford algebra $Cl(p,q)$. The non-commutativity of the new kernel factors constructed from different Clifford algebra square roots of $-1$ allows to place multiple kernel factors to the left and right of the signal function \cite{BHS:History}. The third approach is closely related to the second. It uses a spinor representation for the signal and the square root of $-1$ is chosen as a bivector from the tangent bundle of a higher dimensional immersion of the signal \cite{BBS:colCFT,BB:SpinRepIm}. Except for the work of T. Batard and M. Berthier in ctrb. No. 167 of \cite{KG:ICCA9elproc2011}, we mainly discuss applications related to the second approach.}}

In chp. 8 of \cite{GS:GCompwGA2001} T. B\"{u}low, M. Felsberg, and G. Sommer study harmonic transforms, the quaternionic (Q) Fourier transform (FT), and a type of Clifford FT in $Cl(0,n)$. Furthermore in chp. 9 of \cite{GS:GCompwGA2001} M. Felsberg, T. B\"{u}low, and G. Sommer consider commutative hypercomplex FTs for multidimensional signals, where these commutative hypercomplex algebras can always be regarded as subalgebras of higher dimensional Clifford algebras. 

M. Felsberg et al establish in chp. 10 of \cite{GS:GCompwGA2001} discretizations and fast $n$-D transform algorithms for the QFT of chp. 8 and for the commutative hypercomplex FTs of chp. 9 of \cite{GS:GCompwGA2001}. Finally in chp. 11 of \cite{GS:GCompwGA2001} T. B\"{u}low, and G. Sommer apply local hypercomplex signal representations to image processing and texture segmentation. In \cite{EBC:GCPAS2001} the author applies the QFT to image processing.

M. Bahri, E. Hitzer and S. Adji develop in \cite{BCS:GACEngCS2010}, pp. 93--106, 2D windowed Clifford Fourier transforms, useful for local multivector-image analysis in $Cl(2,0)$. 

In \cite{PGetal:AnVid2dt} a $(x,y,t)$ video is embedded in $Cl(0,3)$. $Cl(0,3)$ can be represented by a pair of quaternions $q_1+\omega q_2$, where $\omega$ with $\omega^2=1$ is the pseudoscalar of $Cl(0,3)$, and $\qi,\qj,\qk$ are the bivectors of $Cl(0,3)$. The B\"{u}low-Felsberg-Sommer Clifford Fourier transform \cite{BFS:ncHypFT} is used to introduce an analytic signal for scalar signals $f(x,y,t)$ in $Cl(0,n)$, and specifically for $n=3$. Then the quaternion pair $q_1+\omega q_2$ is parametrized in terms of a scalar and a pseudoscalar amplitude and 6 phase angles. The analytic signal of $f(x,y,t)=\cos(x) \cos(y) \cos(t)$ is studied and imaged in detail. The authors hope in the future to extract physical information from the phases in medical images.

W. Reich and G. Scheuermann analyze real vector fields by means of Clifford convolution and Clifford Fourier transforms in \cite{BCS:GACEngCS2010}, pp. 121--134. 

R. Bujack, G. Scheuermann and E. Hitzer establish in contribution (ctrb.) No. 141 of \cite{KG:ICCA9elproc2011} a general Clifford algebra framework for most versions of Clifford Fourier transformations in the literature. They show the conditions for linearity and shift theorems and give a range of examples for applying the new general framework.  

T. Ell, who first introduced the quaternion Fourier transform (QFT) for the analysis of 2D linear time-invariant partial differential systems \cite{TE:QFT1993}, gives in ctrb. No. 200 of \cite{KG:ICCA9elproc2011} a current overview of quaternion Fourier transform (QFT) definitions, their relations, inversion, linearity, convolution, correlation and modulation. QFTs have found rich applications in color image processing \cite{SJS:FTcolQuat} by e.g. identifying the colors $r,g,b$ with the coefficients of $\qi,\qj,\qk$.

\cite{HS:ICCA9} and \cite{EH:ICNAAM2011} consider a generalization of the orthogonal 2D plane split (OPS) of quaternions. The OPS split is based on the choice of one or two general pure unit quaternions $\vect{p}^2=\vect{q}^2 = -1$, corresponding to two analysis planes in the $\R^4$ space of quaternions. Vice versa, $\vect{p},\vect{q} \in \HQ$ can be adapted to any set of two orthogonal planes in $\R^4$. This way QFTs can be split into pairs of complex FFTs, QFTs can be interpreted geometrically, new QFTs can be designed (with desired geometric and phase properties) for special applications, and generalizations to higher dimensional Clifford algebras become possible, as e.g. in \cite{EH:QFTgen} to a spacetime FT in $Cl(3,1)$, equivalent to a spacetime multivector wave packet analysis. 

N. Le Bihan and S. Sangwine establish in ctrb. No. 135 of \cite{KG:ICCA9elproc2011} the very first and simplest quaternionic spectral analysis of non-stationary improper (real and imaginary parts may be correlated) complex signals. Based on the 1D QFT a new quaternion valued hyperanalytic signal is established with zero negative frequencies. It is shown how modulus and quaternion angles represent geometric signal features, leading to new notions of angular velocity and complex envelop of complex signals, which are illustrated by examples. 

In ctrb. No. 160 of \cite{KG:ICCA9elproc2011} B. Mawardi generalizes the windowed QFT to a windowed Clifford FT in $Cl(0,n)$, important for local image and signal analysis, and investigates some of its properties. 

In ctrb. No. 167 of \cite{KG:ICCA9elproc2011} T. Batard and M. Berthier introduce a spinor representation of arbitrary surfaces (images). Spinor algebras are natural subalgebras of Clifford algebras (like complex numbers, quaternions, rotor and motor algebras). They investigate applications to image processing, i.e. segmentation, diffusion and a FT for multichannel images to Spin(3). 

In \cite{GSH:RegMImGA} a translation, rotation and scale invariant algorithm for registration of color images and other multichannel data is introduced. In contrast to the original algorithm the proposed algorithm uses the Clifford Fourier transform \cite{ES:CFT} and allows to handle vector valued data in an appropriate way. As a proof of concept the registration results for artificial, as well as for real world data, are discussed.

\subsubsection{Other Clifford integral transformations}

{\change{%
Apart from new types of hypercomplex Fourier transformations, other transformations important for image and signal processing also have higher dimensional generalizations in Clifford algebra. These include amongst others Clifford algebra generalizations of spatially compact atomic function kernels, of the Radon transform, of analytic signal theory, of scale space, of wavelets, of functions on spin groups, of adaptive Fourier decomposition, and of tomographic reconstruction techniques. Several examples are discussed in some more detail in the following.
}}

In \cite{BC:QAFImPro} and ctrb. No. 187 of \cite{KG:ICCA9elproc2011} (by E. Moya-Sanchez and E. Bayro-Corrocha\-no) quaternion atomic function kernels are studied for application in image processing. These kernels are compact in the spatial domain, and they can be adapted to the input signal by broadening or narrowing, ensuring maximum signal-to-noise ratios. They also allow the smooth differentiation of images expressed by quaternion valued functions. The approach also permits to use the generalization of analytic functions to monogenic functions \cite{FS:MonSig}, and to develop a steerable quaternionic multiresolution wavelet scheme for image structure and contour detection. 

Adding a generalized Radon transform allows the detection of color image shape contours. A number of principal experiments with test images, some including color, is performed and the various quaternionic components and phases are discussed. In the future it may be interesting to undertake real world applications and compare them with conventional techniques.  

In ctrb. No. 162 of \cite{KG:ICCA9elproc2011} S. Bernstein generalizes monogenic signal theory (analytic multivector signals) to monogenic curvature scale space, and constructs diffusive wavelets with demodulation applications of 2D AM-FM signals. 

In ctrb. No. 172 of \cite{KG:ICCA9elproc2011} S. Ebert, S. Bernstein and F. Sommen introduce harmonic analysis for Clifford valued functions on the spin group (for example the rotor or motor groups). Secondly, they study Clifford valued diffusive wavelets on the sphere. 

In ctrb. No. 179 of \cite{KG:ICCA9elproc2011} R. Soulard and P. Carre introduce a true color extension of monogenic (i.e. analytic in Clifford analysis) wavelets by embedding color vector signals in the geometric algebra $Cl(5,0)$, different from conventional channel wise (marginal) color image processing. This new multiresolution color geometric analysis with non-separable wavelets yields good orientation analysis well separated from the color information. Statistical coefficient modeling for thresholding and denoising is included. 

In ctrb. No. 197 of \cite{KG:ICCA9elproc2011} T. Qian discusses adaptive Fourier decomposition (AFD), converging fast in energy and point wise. It is a practical solution to the best approximation problem by rational functions of given degree. These results continue to hold for quaternion valued functions. It is used in robust control theory and sound analysis. 

In ctrb. No. 205 of \cite{KG:ICCA9elproc2011} S. Bernstein studies optical coherence tomography (OCD), based on inverse scattering and analytic (generalized to monogenic in Clifford analysis) signals. The usual Green function is modified to include directions. OCT synthesizes a series of adjacent interferometric depth-scans from straight propagation of low-coherence probing beams, decoupling transversal and depth resolutions. And OCT uses backscattering, i.e. light passes twice through the same object.  

\change{%
We finally want to point out the upcoming proceedings volume on quaternion and Clifford Fourier transforms and wavelets \cite{HS:QCFTW}, which has a historic introduction to the development of these transforms (chp. 1) followed by Part I on quaternionic transforms and Part II on Clifford algebra transforms. 
}

\section{Applications to computer and robot vision\label{sc:CompRobV}}

\change{%
Grassmann algebra and Clifford algebra have the advantage that geometric incidence relationships among all geometric objects can universally and exception free be expressed by simple products, instead of often complicated systems of linear equations in the conventional linear algebra approach. This has wide ranging consequences for modeling, design of algorithms, generalizability and computability as well as for speed and accuracy of numerical computations. This applies especially to frequently encountered problems of orientation, pose, motion and tracking, to questions of camera geometry and to the field of scene analysis. For all these areas examples of the application of hypercomplex algebras are discussed in the following three subsections. 
}

\subsection{Orientation, pose, motion and tracking \label{sc:orposmottr}}

Projective geometry has a long history, originated by G. Desargues in the 17th century, and influenced by perspective in art \cite{DN:PengDMath}.  H. Grassmann expressed projective geometry in Grassmann algebra \cite{HG:NBoM1844}. J. Stolfi \cite{JS:OPG} added orientation to projective geometry for computer scientists. D. Hestenes and R. Ziegler \cite{HZ:ProjGeomCA} studied homogeneous $n$D projective geometry in Cl(n+1,0). R. Pappas \cite{ALP:CAcomp1996} introduced orientation to projective geometry in Cl(n+1,0).

The homogeneous model in $Cl(n+1,0)$ of projective space in $\R^{n+1}$, modeling \textit{oriented projective geometry} in $\R^n$ is explained by R. Pappas in \cite{ALP:CAcomp1996} (pp. 233--250), where homogeneous equivalence is restricted to positive scalars in $\R_+$. 

In \cite{SZ:PerAct2000} (pp. 197--207) geometric calculus in $Cl(3,0)$ is employed by H. Sajeewa and J. Lasenby for estimation and tracking of articulated motion described by rotors suitable for forward and inverse kinematics. 

In \cite{SZ:PerAct2000} (pp. 284--293) the even subalgebra $Cl^+(3,0,1)$ (termed motor algebra here) of $Cl(3,0,1)$ is applied to pose estimation by B. Rosenhahn, Y. Zhang and G. Sommer. Moreover, in this motor algebra, an extended Kalman filter design for motion estimation of point and line observations is developed by Y. Zhang, B. Rosenhahn and G. Sommer \cite{SZ:PerAct2000} (pp. 339--348).

In $Cl(3+1,0)$ the projective invariant reconstruction of shape and motion is studied by E. Bayro-Corrochano and V. Banarer in chp. 10 of \cite{BCS:GAappScEng2001}, and by E. Bayro-Corrochano in chp. 7 of \cite{EBC:GCPAS2001}. 

The algebra $Cl^+(3,0,1)$ (also called dual quaternions), where a line is represented by a bivector (equivalent to a quaternion pair) is applied by Y. Zhang, G. Sommer and E. Bayro-Corrochano in chp. 21 of \cite{GS:GCompwGA2001} to motor (motion operator in $Cl^+(3,0,1)$) extended Kalman filtering for dynamic rigid motion estimation from line observations, see also chp. 8 of \cite{EBC:GCPAS2001} by E. Bayro-Corrochano.

In chp. 33 of \cite{DDL:AGACSE2002} by B. Rosenhahn, O. Granert, and G. Sommer, monocular pose estimation of kinematic chains in conformal geometric algebra $Cl(4,1)$ is developed. 

In chp. 37 of \cite{RA:CAappsMPE2004} twist representations of cycloidal curves in $Cl(4,1)$ are used for 3D pose estimation by B. Rosenhahn and G. Sommer.

In part II of \cite{LOS:CAGAapp2005} G. Sommer gives an overview of applications of geometric algebra in robot vision (pp. 258--277): the quaternionic FT, local spectral representations in $Cl(n+1,0)$, and pose estimation of free form objects. G. Sommer, B. Rosenhahn, and C. Perwass give details of the free form object description in $Cl(n+1,1)$ in \cite{LOS:CAGAapp2005}, pp. 278--295, which are also relevant to image processing. R. Wareham, J. Cameron and J. Lasenby describe in \cite{LOS:CAGAapp2005}, pp. 329--349, applications of conformal geometric algebra in computer vision and graphics, including pose and position interpolation, logarithm based interpolation, conics and line images in a para-catadioptric camera. 

D. Gonzalez-Aguirre et al discuss in \cite{BCS:GACEngCS2010}, pp. 299--326, model based visual self-localization using Gaussian spheres in $Cl(4,1)$. 

In \cite{VD:EMfVarGD} a new technique is introduced for estimating conformal geometric algebra $Cl(4,1)$ motion operators (motors) composed of rotators (rotation operators) and translators (translation operators). These operators are all elements of the 8D linear space 
$\mathbb{M}= \text{span}\{1,e_{12},e_{13},e_{23},
e_1n_{\infty},e_2n_{\infty},e_3n_{\infty}, i_3n_{\infty}\}$. The motor manifold $\mathcal{M}$ of rotations and translations is 6D and generated by the 6 bivectors in $\mathbb{M}$. For computing the reciprocal space of $\mathbb{M}$, the vector $n_{\infty}$ in $\mathbb{M}$ needs to be replaced by $n_{0}$, thus the whole conformal geometric algebra $Cl(4,1)$ is needed.  By suitably choosing the conformal object representations of points, spheres, flats (lines and planes), directions, tangents and rounds (circles), the distances and angles (similarity measures) of respective pairs of objects can be represented by scalar products of these objects. 

Symmetrized and weighted sums of these similarities lead to a real scalar cost function, which can be maximized for computing the optimal motor transforming an original set of objects to a corresponding target set of objects. A complete algorithm is developed and principal examples are given. Even in the presence of noise, the optimum Euclidean motion is estimated. Note: (1) This novel algorithm can take geometric data in \textit{any form} (points, $\ldots$, rounds) as input and deliver the motor estimation as output, which makes separate algorithms for each type of object unnecessary. (2) The optimal weighted performance is guaranteed at the same time.

In \cite{FD:RecRBM} a new algorithm is developed for reconstructing rotations and rigid body motions from exact point correspondences through reflections. It can reconstruct orthogonal transformations in Euclidean spaces of any dimension. In non-degenerate situations it can reconstruct orthogonal transformations in spaces with arbitrary metric. A series of reflections aligns the corresponding points one by one. This can be efficiently implemented in conformal geometric algebra, because in CGA the difference of 2 normed conformal points is the bisector plane vector, which can be used as reflection operator, products of these plane vectors give versors for rotations and translations. 

Applications include satellite tracking, and registration of point clouds. A standard matrix algorithm for the determination of \textit{rotations} is almost 2 times slower, quaternions are more than 4 times slower, a Procrustes method was 10 times slower, albeit with the advantage of delivering a least-squares solution in the presence of noise. 

For rigid body motions (rotations and translations) a combination of centroid determination for the translation and the new rotation algorithm turned out to be the fastest, yet the pure versor method was 1.53 times slower than a centroid and rotation matrix method. It works for noisy data, though it may not return a least-squares solution.

In \cite{CL:AttPosTr} a new method for attitude and position tracking of a body based on velocity data is developed, alternative to direction cosine matrices (DCM). The new approach uses (the Lie algebra of) bivectors in geometric algebra as the attitude tracking method of choice, since several features make their performance and flexibility superior to DCMs, Euler angles, or even the classic geometric algebra approach of rotors. Potential advantages of conformal geometric algebra (CGA) are the combination of the integration of angular and linear velocities in a single step, all of which makes bivectors attractive not only for tracking rotations, but for representing general displacements. 

Regarding attitude representations, rotors still require normalization and reprojection, but their integration is computationally efficient, they do not have singularities (of Euler angles), they can be extended to $n$D, and embedded in CGA they even allow solutions for general displacements. All these advantages are shared by the direct bivector approach, yet bivectors only need 3 instead of 4 DOFs, the representation with bivectors is therefore minimal, and it shows excellent optimization by interpolation. Theory and simulation provide compelling reasons to favor the linear space of bivectors over the rotor manifold for tracking attitude from velocity measurements. Especially in its extension to CGA the advantage of bivector update equations is apparent. In addition to computational efficiency, bivectors lead to new integration schemes with dramatically improved performance and have low computational load. 

In \cite{LS:GAfOMCap} an iterative method for motion capture using homogeneous (projective) geometric algebra $Cl(3,0)$ is presented. The external calibration parameters of $n$ cameras observing a scene can be determined (and speedy recalibration after moving the cameras is possible), because in geometric algebra we can easily express geometric entities and differentiate with respect to any element (i.e. multivectors, for example rotors) of the algebra. A byproduct is a fast, efficient and robust world point tracking and reconstruction algorithm after achieving the calibration. A related Bayesian approach is explained in chp. 9 of \cite{EBC:GCPAS2001}. The rotation manifold optimization is very efficient, the parametrization is minimal, there are no singularities (like e.g. Euler angles have), and it can easily be extended to all objects and any dimensions. The results and elements of this algorithm can be easily reused for other algorithms and other scenarios. 

In chp. 8 of \cite{CP:GAAppEng2009} monocular pose estimation in conformal geometric algebra (CGA) $Cl(4,1)$ is presented. Dual quaternions (motor algebra) only allow modeling with lines, in CGA any geometric CGA entity (points, point pairs, lines, circles, planes, spheres, etc.) can be used. 

Other advantages are: No small angle approximation (or exclusion of translations) is needed, due to the linear representation of the Euclidean group in CGA. The constraint equation is quadratic in the motor components and camera model parameters to be estimated. Pose, focal length, and lens distortion of standard lens systems and parabolic catadioptric cameras can be estimated simultaneously, and for both a covariance matrix can be obtained. Thus uncertain data can be used and the uncertainty of the pose can be estimated. 

The optimization algorithms can also be applied to the constraint equations through a tensor representation. The treatment includes a new algorithm for the automatic estimation of the initial pose (no tracking assumptions and no initial pose needed). Alternatively object pose or camera pose can be estimated. Verification experiments were done with a model attached to a robot arm.

\subsection{Camera geometries}

In chp. 7 of \cite{BCS:GAappScEng2001} by E. Bayro-Corrochano and J. Lasenby, and in chp. 4 of \cite{EBC:GCPAS2001} by E. Bayro-Corrochano, the projective (homogeneous) geometric algebras $Cl(2+1,0)$ and $Cl(3+1,0)$ were used for the description of visual geometry of $n$ uncalibrated cameras, greatly simplifying otherwise complex tensorial relations between multiple camera views. Moreover, in chp. 9 of \cite{BCS:GAappScEng2001} a Bayesian inference scheme for camera localization of $n$ cameras and is developed by C. Doran for known range data in $Cl(3,0)$, and for unknown range data in $Cl(3+1,0)$. 

In chp. 14 of \cite{GS:GCompwGA2001} the projective (homogeneous) geometric algebra $Cl(1,3)$ is used by C. Perwass and J. Lasenby for a unified description of multiple view geometry, and further in chp. 15 of \cite{GS:GCompwGA2001} by the same authors, for 3D-reconstruction by vanishing points (points of intersection of parallel lines). In the projective algebras $Cl(3,0)$ for 2D and $Cl(1,3)$ for 3D, intrinsic camera parameters are analyzed and computed based on the absolute conic principle in chp. 16 of \cite{GS:GCompwGA2001} by E. Bayro-Corrochano and B. Rosenhahn. This is applied by the author of \cite{EBC:GCPAS2001} in chp. 7 to visually guided grasping and camera self-localization. Next in chp. 17 of \cite{GS:GCompwGA2001}  coordinate-free $n$D projective geometry for computer vision is developed by H. Li and G. Sommer in $Cl(n+1,0)$. 

In Lecture 5 of \cite{AS:CAapps2004} by J. Selig, biquaternions (dual quaternions), isomorphic to $Cl^+(0,3,1)$, are used for a two camera correspondence problem. 

In chp. 7 of \cite{CP:GAAppEng2009} C. Perwass  discusses an inversion camera model, followed by single camera monocular pose estimation in chp. 8. 

T. Debaecker, R. Benosman and S. Ieng study in \cite{BCS:GACEngCS2010}, pp. 277--298, an image sensor model using conformal geometric algebra $Cl(4,1)$ for cone-pixel camera calibration and motion estimation. 

In \cite{JLetal:NGeoMfCV} a technique for estimation of camera (or object) motion and scene structure from 2 scene projections, range data known (or unknown) is presented, using the advantages of geometric algebra $Cl(3,0)$, already described in the discussion of \cite{LS:GAfOMCap} in section \ref{sc:orposmottr}. Motion is analyzed from a pair of images, even for the case, where only 2D information is available without range information in 2D projective geometric algebra $Cl(3,0)$. For unknown range data both motion and depth coordinates of image points in the presence of measurement uncertainty are computed iteratively. 

The formulation directly provides least-squares estimates of the motion and structure simultaneously, with analytic multivector derivatives for the whole set of unknowns in the problem. Note that the direct minimization with respect to rotors is not readily available in frameworks other than geometric algebra. The results are verified by simulations in comparison with basic linear algorithms for structure and motion estimation. The accuracy is improved, in some cases by large factors.

\subsection{Scene analysis \label{sc:ScnAnalysis}}

In chp. 29 of \cite{CgeoAPME1996} (by J. Lasenby) the homogeneous model in $Cl(1,3)$ is used for projective geometry in $\R^{1,3}$ modeling 3D space to describe geometric invariance in computer vision. 

H. Li et al in \cite{LOS:CAGAapp2005}, pp. 363--382 develop a system for $n$D object representation and detection from single 2D line drawings in (projective) Grassmann  algebra. H. Li, L. Zhao, and Y. Chen in \cite{LOS:CAGAapp2005}, pp. 383--402, use Grassmann-Cayley algebra for polyhedral scene analysis (realizability of 2D line drawings of 3D polyhedra) combining parametric propagation with calotte analysis. Y. Wu and Z. Hu develop with projective geometric invariants in \cite{LOS:CAGAapp2005}, pp. 403--417, a unified and complete framework of projective invariance for six points, important for vision tasks and critical geometric information, reducing the need for explicit estimation of the camera projective matrix and the optical center. 

C. Perwass discusses in chp. 6 of \cite{CP:GAAppEng2009} the construction and estimation of uncertain geometric entities and operators in geometric algebra. 

A new Grassmann algebra (restriction of projective Clifford algebra to the outer product) based framework for global $n$D visibility computations is developed in \cite{ACFM:FnDVis}. 
Given a set of obstacles in the Euclidean space, two points in the space are said to be visible to each other, if the line segment that joins them does not intersect any obstacles \cite{Wiki:Visibility}.
All rendering algorithms aim at simulating the light transfer in a virtual environment, which strongly depends on the mutual visibility of each element in the scene. 
This is fundamental in computer graphics, well understood in 2D, but so far not in $n$D. 

Grassmann algebra allows a high level of abstraction for visibility problems, general applicability, and no exceptions, thereby superseding Stolfi's line geometric representation \cite{JS:OPG}. In Grassmann algebra the exterior (outer) product achieves line classification, equivalent to point versus hyperplane classification -- well-defined and computationally robust. Lines stabbing an $n$-dimensional convex face are characterized. The stabbing lines set is the intersection of the line bivector blade set and a convex polytope. The latter is minimal, leading to algorithmic improvements. 

The new framework is illustrated with computations of soft shadows for 3D illuminated scenes and compares favorably with conventional production rendering software, as can be clearly seen in soft shadow details. The new visibility framework can be considered as a black box and easily plugged into any applications that need to perform visibility queries.

\section{Applications to kinematics and dynamics of robots\label{sc:KinDynRob}}

\change{%
Due to the pioneering work of D. Hestenes \cite{DH:STA,DH:NFCM} Clifford's geometric algebra has been recognized as ideally suited for mechanics, including the fields of kinematics and dynamics. This has been met with great interest in the robotics community \cite{DFM:GAfCS}. Geometric algebra allows to formulate algorithms on a very high geometrically intuitive level, which additionally leads to very efficient software and hardware implementations \cite{HLSK:EffIK_CGA}. The following subsections on kinematics (and dynamics) feature a number of concrete examples, where geometric algebra served as the principal modeling language, leading to new algorithms with efficient implementations. 
}

\subsection{Kinematics}

In chp. 29 of \cite{CgeoAPME1996} geometric algebra $Cl(3,0)$ is extended by J. Lasenby to hyperspinors \cite{DH:HypSpin} by including a new element $\epsilon$, with $\epsilon^2=0$, similar to dual quaternions (one form of biquaternions). This allows to linearize translations like rotations and create the unified notion of motion operator (motor) with application to a robot arm of $n$ links. 

In \cite{SZ:PerAct2000} (pp. 22--47) the homogeneous model in $Cl(n+1,0)$ is used by L. Dorst and R. van den Boomgard to develop a system theory of contact important for robot interactions. This approach is further extended by L. Dorst in chp. 17 of \cite{BCS:GAappScEng2001} to $Cl(n+1,1)$. 
In \cite{SZ:PerAct2000} (pp. 115--133) a Lie model for $n$D Euclidean geometry in $Cl(n+1,2)$ is developed by H. Li for geometric contact problems of spheres and hyperplanes. It is a supermodel for the conformal model in $Cl(n+1,1)$, which is here also called homogeneous model, since all its points and blades are also homogeneous. 

In \cite{SZ:PerAct2000} (pp. 104--114) motion is modeled by J. Lasenby, S. Gamage and M. Ringer in $Cl(3,0)$ for tracking, analysis and inverse kinematics, and information is given, that for complex geometries conformal geometric algebra in $Cl(4,1)$ will be employed. 

In chp. 11 of \cite{BCS:GAappScEng2001} the Clifford algebra $Cl(0,3,1)$, which is related to Clifford's biquaternions $Cl^+(0,3,1)$ (dual quaternions), is used by J. Selig for inverse kinematics of a serial manipulator. 
Moreover, in Lecture 5 of \cite{AS:CAapps2004} biquaternions , isomorphic to $Cl^+(0,3,1)$ are applied by J. Selig to the inverse kinematics of a robot arm. 

In chp. 12 of \cite{BCS:GAappScEng2001} a robot task trajectory design algorithm with Bezier interpolation in $Cl^+(0,4)$, here also called double quaternions, is developed by S. Ahlers and J. McCarthy. 

Furthermore, bivector Lie algebra embedded in conformal geometric algebra $Cl(n+1,1)$ and algebraic incidence relations are studied by E. Bayro-Corrochano and G. Sobczyk in chp. 13 of \cite{BCS:GAappScEng2001}, by E. Bayro-Corrochano in chp. 3 of \cite{EBC:GCPAS2001}, and by E. Bayro-Corrochano, P. Lounesto, and L. Lozano in chp. 32 of \cite{DDL:AGACSE2002} with applications to robotics and image analysis. 

In chp. 18 of \cite{GS:GCompwGA2001} the geometry and algebra of robot kinematics is developed by E. Bayro-Corrochano in $Cl(3,0)$, $Cl(3,0,1)$ and the motor algebra $Cl^+(3,0,1)$, compare also \cite{EBC:GCPAS2001}. The approach in the motor algebra $Cl^+(3,0,1)$ is then applied by E. Bayro-Corrochano and D. K\"{a}hler in chp. 19 of \cite{GS:GCompwGA2001} to direct and inverse kinematics of a robot manipulator, see also \cite{EBC:GCPAS2001}. Moreover the algebra $Cl^+(3,0,1)$ (also called dual quaternions), where a line is represented by a bivector (equivalent to a quaternion pair), is used by K. Daniilidis in chp. 20 of \cite{GS:GCompwGA2001} for motion alignment. 

D. Hestenes in \cite{BCS:GACEngCS2010}, pp. 3--34, develops new tools for computational geometry and screw theory in conformal geometric algebra $Cl(4,1)$ with applications to linked rigid bodies and robotics. 

D. Hildenbrand, J. Pitt and A. Koch describe in \cite{BCS:GACEngCS2010}, pp. 477--494, an application of high performance parallel computing based on conformal geometric algebra with Gaalop to the inverse kinematics of the leg of a humanoid robot. 

A novel forward and backward reaching inverse kinematics (FABRIK) solver implemented in conformal geometric algebra (CGA) is developed in \cite{AL:IKSuCGA}. It can find joint positions by locating points on lines. It can treat most joint types and support biomechanical constraints on chains with single and multiple end effectors. The example of a human hand is employed for efficient pose tracking and reconstruction in real time, with smooth and natural motion results. Data from a markered optical motion capture system are used. In the forward reaching phase, an end effector is moved to a target position.  Then length of joint radius spheres, and new end effector point to old opposite end of joint lines, are intersected to construct the new opposite end point of each joint. This is done for every joint in a chain, finally relocating the root. In the 2nd backward reaching phase the root point is moved to its initial position and the algorithm repeated in reverse sense. 

The approach is suitable for computer animations and provides a framework for many IK applications, computer vision and robotics. Rotors (rotation operators in CGA) are found more numerically stable and more efficient than rotation matrices. CGA can model both virtual and mechanical objects. The human hand model used has a highly constrained feasibility set, with both rotational and orientational constraints. For example, joint points on palm plane constraints lead to direct estimates of joint positions by intersecting inter-joint distance spheres and the palm plane. The algorithm is geometrically intuitive and compact. 

Experiments with a 10 camera phase space motion capture system, data capture at 100 Hz, were carried out. The CGA method processed up to 70 frames per second in MATLAB. In only 1.43 ms per frame 25 joints were fitted accurately with FABRIK based on 8 markers, ensuring normal movements without asymmetry, irregular bends and rotations. No oscillations or discontinuities occurred. 

The inverse kinematics of a robot arm is intuitively formulated in conformal geometric algebra (CGA) $Cl(4,1)$ in \cite{DH:GCinCG} in four simple steps. The robot arm consists of three links, three joints and a gripper with a total of five joint angle degrees of freedom. Every step simply uses basic elements of CGA: points, lines, planes and spheres, their combinations and intersections, which are elementary algebraic products in CGA. Finally inner products of the thus constructed planes and lines result in the joint angles. Every step can be fully visualized with the free CLUCalc \cite{CP:CLUCalc}. The paper also shows how to best fit planes and spheres in a set of points.

\subsection{Dynamics}

Many results of robot kinematics in the preceding section also apply to robot dynamics in Clifford algebra. 

Spacetime algebra (STA, \cite{DH:STA}) in $Cl(1,3)$ is used by S. Gull, C. Doran and A. Lasenby in chp. 7 of \cite{CgeoAPME1996} to treat the dynamics of rigid bodies and of elastic media. 

Homogeneous 3D rigid body mechanics with elastic coupling is studied by D. Hestenes in conformal geometric algebra $Cl(4,1)$ in chp. 19 of \cite{DDL:AGACSE2002}.

Moreover in $Cl(3,0)$ dynamical information is inferred from 3D position data by H. Sajeewa and J. Lasenby in chp. 35 of \cite{DDL:AGACSE2002}.

\section{Applications to control problems\label{sc:CtrlProbl}}

\change{%
In the development of applications of geometric algebra it has been found parti\-cularly useful, that one can differentiate with respect to versors. In the conformal model versors include rotors, translators and most general motion operators (motors). Many objects (vectors, points, lines, spheres, motors) can all be represented in one common algebra, simplifying software implementations. Linear operators can be universally extended via outermorphisms, leading to compact solutions. Both problem formulation and solution become more elegant than in equivalent matrix theoretic formulations \cite{VN:CalTPosCGA}. In the following two subsections application examples of geometric algebra to platform control, optical positioning and geometric constraint solving are considered.
}

In chp. 36 of \cite{DDL:AGACSE2002} the singularities of inline planar platforms are studied by M. Baswell, R. Ab{\l}amowicz and J. Anderson in the algebra $Cl^+(3,0,1)$. In Lecture 5 of \cite{AS:CAapps2004} biquaternions (dual quaternions), isomorphic to $Cl^+(0,3,1)$ are used by J. Selig for describing a Stewart platform. 

\subsection{Optical positioning system}

In \cite{VN:CalTPosCGA} a new algorithm is developed for the rapid calibration of 3D positions of multiple stationary point targets, part of an optical positioning system using \textit{conformal geometric algebra} (CGA). An earlier algorithm in the homogeneous (projective) model \cite{VLK:SelfCalGA} is thus extended to CGA. CGA made expressions simpler, simplified the software implementation because it represents many objects and operators (vectors, points, lines, dual spheres, motors) in a single algebra, and the use of the outermorphism concept (related to the universal property of Clifford algebras) allowed a very compact representation of the results. The representation of inverses was simplified (over matrix theory) by using duality, and inverses could be more efficiently differentiated, e.g., for computing the closest point to a number of lines. 

In the framework of CGA a nonlinear solution is developed. This includes a coordinate free approach to differentiating rotors applied to motion rotor differentiation. Separate differentiation with respect to a rotor and a translation vector became thus unnecessary. The coordinate-freeness allowed to delay and change the representation choice easily. 
 
In the experiment a group of rigidly co-located calibrated cameras is moved to several positions, and target images are acquired. Target pixel coordinates are extracted and transformed into 3D camera-target lines, and then used as algorithm input data. Experimental results show good robust performance even in the presence of noise. For better noise (e.g. due to subpixel target locations) modeling, the current simple and elegant error function should be further improved. 

\subsection{Geometric constraint solving}

In \cite{SAY:CGAGeoConstr} the use of conformal geometric algebra (CGA) in geometric
constraint solving (GCS) is studied for declarative modeling (DM) for mechanical systems and mechatronics in parametric computer aided design (CAD). DM explicitly describes properties of designed objects (primitives of points, lines and planes, etc.), but not how to construct them, a constraint solver is used for that. In DM the description step is followed by the generation step, which explores the solution space. The final lookup step visualizes the results. CGA allows to algebraically specify geometric objects, their mutual constraints (design requirements of size, position and orientation) and constraint solving. CGA can thereby postpone the introduction of coordinates until concrete visualizations are needed. 

Thus in CGA a new \textit{unified framework} for topologically and technologically related surfaces (TTRS), their underlying Lie algebras, and their screw-based representations can be constructed, suitable for symbolic GCS and algebraic classification. Moreover CGA offers a new language for chirality  and mobility specifications, because e.g. rotators and geometric constraints are both described in CGA. The scientific universality of geometric algebra in physics allows to incorporate physical interactions with geometry in DM for mechatronics. ~\cite{SAY:CGAGeoConstr} includes concrete CGA application examples for symbolic GCS and for geometric classification.

\section{Other applications\label{sc:othApp}}

\change{%
In order to additionally demonstrate that Clifford's geometric algebra can revolutionize and advance the fundamental knowledge in mathematics, quantum physics and software development we review in this section unimaginable progress in automatic but human readable geometric theorem proving, applications of geometric algebra to quantum computing, and finally briefly survey Clifford algebra software. 
}

{
\subsection{Automatic geometric theorem proving}

In \cite{HL:IAandGR} automatic geometric theorem proving in Grassmann algebra, Grassmann-Cayley algebra, inner product Grassmann algebra and Clifford algebra  is discussed. Aiming for a complete system of advanced orthogonal invariants naturally leads to Clifford algebra. Using conformal geometric algebra provides a unified framework including projective, affine, Euclidean, hyperbolic, spherical, elliptic and conformal geometries. Geometric theorems in these geometries can now be proved automatically by computers resulting in short (several lines) human readable proofs, with complete control of middle expression swell, so that most expressions only involve monomials and binomials. A host of explicit examples is given. This approach is also of great practical value in subjects like scene recognition, compare section \ref{sc:ScnAnalysis} where H. Li et al's contributions in \cite{LOS:CAGAapp2005} are discussed. 

\subsection{Quantum computing and light}

Quantum computing in geometric algebra is discussed in \cite{DDL:AGACSE2002}, regarding one, two (chp. 20 by R. Parker and C. Doran) and multiparticle quantum systems (chp. 21 by T. Havel, C. Doran) in multiparticle spacetime algebra $Cl(n,3n)$, $n=$ number of particles, see also \cite{RA:CAappsMPE2004}. The treatment includes entanglement (chp. 21), and a discussion of the brain as a $Cl(p,q,r)$ Clifford algebra quantum computer (chp. 25 of \cite{DDL:AGACSE2002} by V. Labunets, E. Rundblad and J. Astola). 

Light polarization (chp. 26 of \cite{DDL:AGACSE2002} by Q. Sugon and D. McNamara) is discussed in spacetime algebra $Cl(1,3)$. In \cite{RA:CAappsMPE2004} (chp. 30 by A. Lasenby, C. Doran and E. Arcaute) a unique algorithm is developed for \textit{real time} computations of electromagnetic fields scattered from conductors of arbitrary shape, a truly first principle method for computing fields in electromagnetic cavities and lighting on all scales. 

\subsection{Clifford algebra software}

Regarding software, in \cite{ALP:CAcomp1996} the editors explain in the preface, that since 1971 the following softwares have been used for numeric and symbolic Clifford algebra computations: FORTRAN, C++, PASCAL, CLICAL, muSIMP, LISP, muMATH, MCSYMA, MAPLE, MATHEMATICA, DERIVE, AXIOM, SCRATCHPAD, MATLAB, and REDUCE. Several of these implementations are discussed in some detail. Moreover the editors of \cite {AS:CAapps2004} briefly review Clifford algebra software (up to 2004) in Appendix 7.1. 

In \cite{RA:CAappsMPE2004}, chp. 35, the implementation of a Clifford algebra $Cl(p,q,r)$ co-processor design on a field programmable gate array (FPGA) is explained and evaluated up to $n=p+q \leq 8$, by C. Perwass, C. Gebken and G. Sommer. C. Perwass discusses in chp. 5 of \cite{CP:GAAppEng2009} the numerical questions of geometric algebra. A. Rockwood and D. Hildenbrand discuss in \cite{BCS:GACEngCS2010}, pp. 53--70, the computational efficiency of geometric algebra using the Gaalop environment. D. Hildenbrand, J. Pitt and A. Koch describe in \cite{BCS:GACEngCS2010}, pp. 477--494, in further detail how the high performance parallel computing based on conformal geometric algebra with Gaalop is designed.  
}

%We finally want to direct the interested reader to a new book \cite{DH:FGAComp} by D. Hildenbrand on Foundations of Geometric Algebra Computing\footnote{\change{We thank one of the anonymous reviewers for the information about this upcoming new book.}}, which focuses in Part I on the mathematical foundations; whereas Part II explains the interactive handling of geometric algebra; and Part III deals with computing technology for high-performance implementations based on geometric algebra as a domain-specific language in standard programming languages such as C++ and OpenCL. 

\section{Conclusions}

We have given an overview on how applications of Clifford's geometric algebra in neural computing, image and signal processing, computer and robot vision, control problems and other areas have developed over the past 15 years. Apart from special applications (like visibility, where Grassmann algebra is sufficient) nowadays conformal geometric algebra $Cl(4,1)$, or more general $Cl(p+1,q+1)$, has become the most widely applied system. Complex numbers, quaternions, (complex and dual) biquaternions and motor algebra are all included as subsystems (subalgebras) in conformal geometric algebra. 

In many cases better accuracy, higher speed, and a wider exception free and singularity free scope more than justify the use of geometric algebra, apart from the great conceptual gains of intuitive algorithmic methods. 

\section*{Acknowledgments}

E.H. gratefully acknowledges God: Soli Deo Gloria (J.S. Bach). He further thanks his family for gracefully tolerating his work even during holidays and vacations. He thanks J. Koenderink for a helpful comment. We warmly thank the anonymous reviewers for several corrections and helpful comments.

\end{document}